# A Versatile Queuing System For Sharing Economy Platform Operations

SONG-KYOO KIM & CHAN YEOB YEUN

## ABSTRACT

The paper deals with a sharing economy system with various management factors by using a bulk input G/M/1 type queuing model. The effective management of operating costs is vital for controlling the sharing economy platform and this research builds the theoretical background to understand the sharing economy business model. Analytically, the techniques include a classical Markov process of the single channel queueing system, semi-Markov process and semi-regenerative process. It uses the stochastic congruent properties to find the probability distribution of the number of contractors in the sharing economy platform. The obtained explicit formulas demonstrate the usage of functional for the main stochastic characteristics including sharing expenses due to over contracted resources and optimization of their objective function.

## 1. INTRODUCTION

This paper deals with a sharing economy platform with critical operational factors which include a status of the owners (contractors, suppliers) and seekers (subscribers, customers) under the sharing economy. Basically, the sharing economy platform makes the collaborative consumption by the activities of sharing, exchanging, and rental of resources without owning the goods [5]. A growing concern about climate change have made the collaborative consumption as appealing alternative for consumers [11]. This concept becomes a broader sense to describe any sales transactions which are done in online market places, even in business-to-business (B2B). This sharing economy began to spread widely by sharing the unused resources between individuals. The startup companies including Uber and AirBnB are not only providing technological platforms to facilitate transactions but also become real-world companies, with the same responsibilities as transportation companies under the sharing economy platforms [15][20]. The constituent factors of a sharing economy business model could be largely divided into eight factors for the business operations: value proposition; financial profit and loss; resources; process; target customers; exterior cooperation; and logistical flow [7]. The operations must evolve along with the new sharing business model and they should preserve the value created through the innovation. The paper deals with modeling the business operations between resource owners (suppliers) and consumers (customers). As a platform provider, a manager should control the balance between them. The paper mainly suggests that resource owners have mainly two types of the contract in the sharing economy platform. Individual contractors are typical suppliers in the system which are independent and have less obligations. But company contractors who are directly reserved by a platform provider have more obligations to support goods. The relationship between a platform provider and various types of suppliers in the sharing economy system is described in Figure 1.

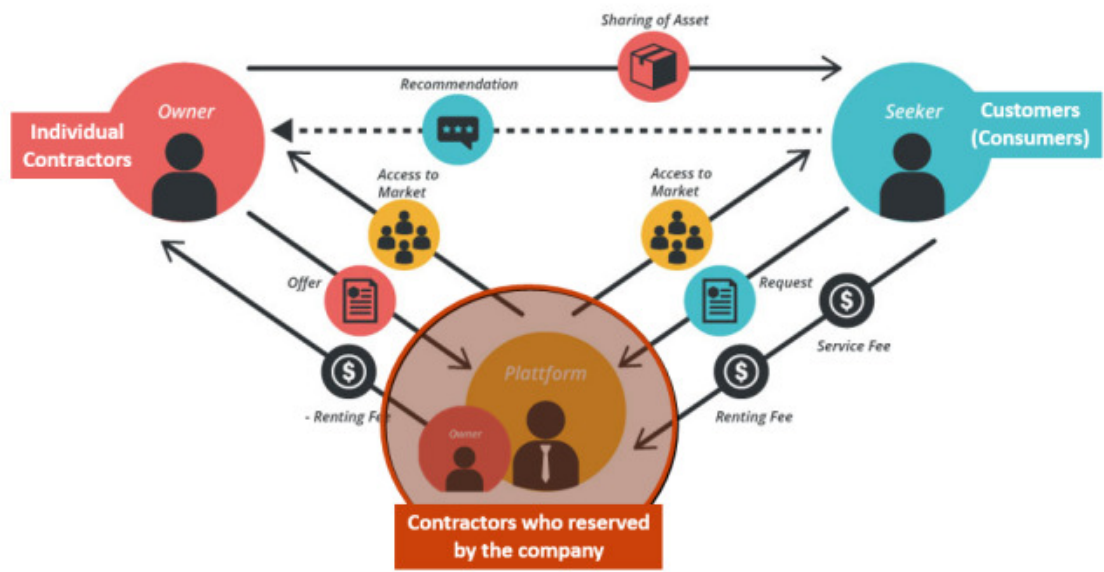

**Figure 1.** Sharing Economy with Two Types of Contractors [4]

This type of a sharing economy system could be described as a classical G/M/1 type queueing model and its operations are critical for all companies in each section of the sharing economy. The main target of the sharing economy operations is balancing conflicted goals such as keeping enough for various types of supplying resources (e.g., car drivers in Uber, room owners in AirBnB). The literature reviews of queueing research have contained massive studies regading queues with various situations [2][8][14][16] and the $M^b/G^b/1$ type queues have been not only theoretically studied [3][22] but also studied based on the algorithmic approach [1].

The model in the paper considers two types of contracts for suppliers. The individual contractors is typical resource owners in a sharing economy system. As it mentioned, each contract is individual which is similar with a part-time contract. The other contract is rather permanent than first one and the contractors will be directly controlled by the service providers. In the terms of operating cost, the individual contract is cheaper but it has the risk for out of the service at the same time. Therefore, the contract type of resource owners should be balanced. Notably, the valid suppliers are posted on the (online) platform as a bulk of suppliers after gathering and screening. The duration of gathering (and screening) resource owners is considered as a general inter-arrival time in atypical single channel queue and the service which is consumed by customers is based on the Markovian service time. This sharing economy model is well described as a $G/M/1$ type queueing model rather than $M/G/1$ type queue [17]. Hence, the paper deals a class of $G^x/M/1$ queueing model with the continuous time perspectives. More specifically, a $G/M/1$ type system with a finite capacity $w$ is considered and each input (i.e., offering resources from the individual contractors) is stored as a fixed batch size $v$. A basic $G^v/M/1$ type queue has been well studied [5-6] in various area including inventory controls [9][16] and manufacturing systems [13] but it has never been adapted for sharing economy models. The "congruent properties" explored in the paper gives an explicit relationship between two models in terms of probabilities for the queueing processes. Explicit obtained formulas demonstrate a relatively effortless usage of functionals for the main stochastic characteristics and their objective function supports the optimization of the system. The main objective of this article is to analyze the random process describing the sharing economy system at any moment of time in equilibrium. The versatile scientific approaches by using modern technologies including an artificial intelligence and a complex information system have been adapted for resolving social and economic challenges [10][18]. The theoretical approach to find the optimization of a sharing economy platform might be a basic clue for analyzing social and economic challenges to understand vital factors of the problems before adapting complex intelligent systems.

The article is organized as follows: Section 2 presents a "congruent model" which has the same characteristics as the principal model. Some expressions provide a very simple connection between the main and the congruent models. This section also includes mathematical formalism of two models and the congruent properties are further extended to worthwhile relationship between Model 1 ($G^v/M/1/w$) and Model 2 ($M/G^v/1/w$). Section 3 deals with the Markov processes and continuous time processes for Model 2 (which bring back to Model 1 via the stochastic congruence). In Section 4, functions of the semi-regenerative process are introduced to demonstrate tractability of the main results, followed by optimization of a relevant objective function. Section 5 deals with optimization of atypical sharing economy platform operation case which concluded with numerical examples.

## 2. STOCHASTICALLY CONGRUENT MODELS

The stochastically congruent model which referred as Model 2 is introduced which is simpler than the main model (i.e., Model 1). The stochastically congruence is widely applied to solve a problem by using a well known model which is connected with the original problem. In this paper, Model 2 is similar to Model 1 and the analytical solution of Model 2 are already well known [6][19]. According to the congruency, the states of Model 1 are equivalent to states of Model 2 which are simply flipped over and counted other way around. In other words, Model 2 is directly connected with Model 1. All above models (Model 1 and 2) could be formally described.

Denote by $Z_t^1$ the total number of contracts from resource suppliers $w$ at time $t$ in Model 1. Once one batch of individual contract $v$ from a supplier is ready, the resource is published for customers. As it is mentioned, this sharing economy model could be described as $G^v/M/1/w$ queueing model. Let $\tau_0 (=0), \tau_1, \ldots$ be the moment of a batch contract is published on the system (the online platform) and let $D_1, D_2, \ldots$ be the duration of publishing bulk service suppliers. The random variables $\{D_n\}$ are iid (independent and identically distributed) with a common CDF (Cumulative Density Function)

$$A(x) := P\{D_n \leq x\},\ x \geq 0, \tag{2.1}$$

with the mean $a = \mathbb{E}[D_n] < \infty$. Each of consuming services (or users get the service from the contracted suppliers) is independent each other and exponentially distributed with the parameter $0 < \lambda < \infty$. The summary of the related values could be provided by the table (see Table 1).

| Variables | Descriptions in the sharing economy platform |
|---|---|
| $v$ | The number of the individual contractors as a service supplier |
| $w$ | The total number of contractors including company reserved |
| $s = (w-v)^+$ | The number of the company reserved contractors |
| $a = \mathbb{E}[D_k]$ | Duration between posting the suppliers on the platform |
| $\lambda$ | Customer arrivals to use the services on the platform |

**Table 1.** Notation Summary

Model 2 describes the number of customers in a $M/G^v/1/w$ queueing system which is a single sever process with the fixed batch size $v$. It is like a single channel queue $M/G^v/1/\infty$ except for the finite waiting room $w$ [6]. Let $T_0(=0), T_1, \ldots$ be the moment of contracted suppliers are published on the platform and it is assumed that the last action does not affect the status of the system in this particular problem setting. The random variables $T_{n+1} - T_n$ (in Model 2) are stochastically equivalent to $\tau_{n+1} - \tau_n$ of Model 1. It can be shown that $\tau_1, \tau_2 \ldots$ and $T_1, T_2, \ldots$ are sequences of stopping times of the processes

$$\left(\Omega_1, \mathfrak{U}_1, (P^x)_{x \in E_1}, (Z_t^1; t \geq 0)\right) \to E_1 = \{0, 1, \ldots, s\} \tag{2.2}$$

and

$$\left(\Omega, \mathfrak{U}, (P^x)_{x \in E}, (Z_t; t \geq 0)\right) \to E = \{0, 1, \ldots, s\}, \tag{2.3}$$

respectively, and that these processes are regenerative relative to these sequences. Let us assign:

$$\xi_n := Z^1_{\tau_n-} \text{ and } X_n := Z_{T_n-},\ n = 0, 1, 2, \ldots.$$

Consequently,

$$\left(\Omega_1, \mathfrak{U}_1, (P^x)_{x \in E}, (\xi_n; n = 1, 2, \ldots)\right) \to E \tag{2.4}$$

and

$$\left(\Omega, \mathfrak{U}, (P^x)_{x \in E}, (X_n; n = 1, 2, \ldots)\right) \to E \tag{2.5}$$

are embedded Markov Chains (MCs). Both Markov chains are stochastically congruent. Their limiting probabilities are expressed through the common invariant probability measure $\mathbf{P} = (P_0, P_1, \ldots, P_w)$ where $P_k = P\{X_n = k\}$ and $\mathbf{P^1} = (P_0^1, P_1^1, \ldots, P_w^1)$ where $P_k^1 = P\{\xi_n = k\}$. The states of the $G^v/M/1/w$ queue are equivalent to the states of the flipped $M/G^v/1/w$ queue. Therefore, the state $k$ in $G^v/M/1/w$ queueing system "corresponds" to the state $w-k$ in $M/G^v/1/w$ queueing system which yields that

$$P_k^1 = P_{w-k}, \; k = 0, \ldots, w. \tag{2.6}$$

Since $(\xi_n)$ and $(X_n)$ are stochastically congruent, only one of them will be treated in the next section.

## 3. Embedded Processes and Continuous Parameters of Model 2

The embedded process $\{X_n\}$ of Model 2 is a time-homogeneous Markov chain and $X_n$ is the number of customers at time $(T_n-)$ and it has two types of transition matrices which depend on the relationship between the batch size $v$ and the system capacity $w$. The transition probabilities could be changed because of these two factors but the method of solving stationary probabilities are essentially same. Let $V_{n+1}(j)$ be the number of completely processed customers during the period $[T_n, T_{n+1})$ when $X_n = j - 1$. Then we have:

$$X_{n+1} = \begin{cases} (X_n - x)^+ + V_n, & X_n = 0, 1, \ldots, v-1, \\ V_n, & X_n = v, v+1, \ldots, w \end{cases} \tag{3.1}$$

Two types of the transition probabilities $p_{jk} := P\{X_{n+1} = k | X_n = j\}$, $j = 0, 1, \ldots, w$ and $k = 0, 1, \ldots, w$, are as follows, if $v \leq w - v$ (Type-1):

$$p_{jk} = \begin{cases} \psi(k), & 0 \leq j \leq v, \\ \psi(k - j - x), & v < j \leq w, j \leq k \leq w - v, \\ 0, & \text{otherwise}, \end{cases} \tag{3.2}$$

or, if $v > w - v$ (Type-2):

$$p_{jk} = \begin{cases} \psi(k), & 0 \leq j \leq v, 0 \leq k \leq w - v \\ 0, & \text{otherwise}. \end{cases} \tag{3.2a}$$

where

$$\psi(k) = \int_{\mathbb{R}_+} e^{-\lambda x} \cdot \frac{(\lambda)^k}{k!} dA(x), \tag{3.3}$$

Both TPMs (Transition Probability Matrices) of Model 2 (Type 1 and 2) repent a finite sub-matrix which occurs within the corresponding TPM of $M/G^v/1/\infty$ queueing system. The embedded probabilities $\mathbf{Q} = (Q_0, Q_1, \ldots)$ of a $M/G^v/1/\infty$ queue are as follows [5]:

$$Q(z) = \sum_{i=0}^{\infty} Q_i z^i = \frac{\alpha(\lambda - \lambda z)\left[\sum_{i=0}^{v-1} Q_i(z^v - z^i)\right]}{z^b - \alpha(\lambda - \lambda z)}, \tag{3.4}$$

where

$$\alpha(\theta) = \mathbb{E}\left[e^{-\theta D_k}\right] = \int_{\mathbb{R}_+} e^{-\theta u} dA(u). \tag{3.5}$$

The generating function $Q(z)$ converges inside the unit open disc centered at zero. Alternatively, $\mathbf{Q}$ is also the solution of equations $\mathbf{Q} = \mathbf{Q} \cdot \mathbf{M}_\infty$ and $(\mathbf{Q}, \mathbf{1}) = 1$ where $\mathbf{M}_\infty$ is the TPM of the infinite queueing system. It is noted that $\mathbf{M}_{(w-v+1)}$ is the $(w+1) \times (w+1)$ finite sub-matrix of $\mathbf{M}_\infty$. Let the vector $\widehat{\mathbf{Q}} = \left(\widehat{Q}_0, \widehat{Q}_1, \ldots, \widehat{Q}_{w-v}, 0, \ldots, 0 \, (= \widehat{Q}_w)\right)$ be the solution of the equation $\widehat{\mathbf{Q}} = \widehat{\mathbf{Q}} \cdot \mathbf{M}_{(w+1)}$. Because of the linearity of the equations $\mathbf{Q} = \mathbf{Q} \cdot \mathbf{M}_\infty$ and $\widehat{\mathbf{Q}} = \widehat{\mathbf{Q}} \cdot \mathbf{M}_{(w+1)}$, $\widehat{Q}_0, \widehat{Q}_1, \ldots, \widehat{Q}_{w-v}$ differ from $Q_0, Q_1, \ldots, Q_{w-v}$ by the same constant respectively. It means a constant $A$ (proportionality) exists such that

$$A = \frac{\widehat{Q}_0}{Q_0} = \frac{\widehat{Q}_1}{Q_1} = \cdots = \frac{\widehat{Q}_{w-v}}{Q_{w-v}}. \tag{3.6}$$

Since we will use the probabilities $\widehat{\mathbf{Q}} = (\widehat{Q}_0, \widehat{Q}_1, \ldots, \widehat{Q}_{w-v})$ for the solutions of $M/G^v/1/w$, the constant $A$ ($> 0$) will be determined from $(\mathbf{P}, \mathbf{1}) = 1$, where $\mathbf{P}$ is the invariant probability measure of the embedded process in $M/G^v/1/w$ system. The probabilities $P_k$, $k = 0, 1, \ldots, w$ were obtained [21] from the expansion of the generating function (3.4). From (2.1) and (3.1)-(3.6), it follows that:

$$P_k = \begin{cases} A \cdot Q_k, & k = 1, \ldots, w - v, \\ 0, & k = w - v + 1, \ldots, w, \end{cases} \tag{3.7}$$

where

$$A^{-1} = \sum_{i=0}^{w-v} Q_i (\leq 1), \tag{3.7a}$$

and

$$P(z) = \sum_{k=0}^w P_k z^k = \sum_{k=0}^w \sum_{j=0}^w p_{jk} P_j z^k. \tag{3.8}$$

The next step is finding the continuous time parameter queueing process of Model 2. The below treatment is similar to that of Kim and Dshalalow [12]. Let $(N_t)$ be the counting process associated with the point process $(T_n)$ and $V_t := T_{N_t+1} - t$, $t \geq 0$. $\{V_t; t \geq 0\}$ gives the residual time from time $t$ to the next service completion $T_{N_t+1}$. The process

$$(\Omega, \mathfrak{F}, (P^x)_{x \in E}, (Z_t, V_t)_{t \geq 0}) \to (\Phi, \sigma(\Phi)), \Phi := E \times \mathbb{R}_+,$$

is weak Markov. Let $(P_t^i)$ be the Markov semi-group associated with $(Z_t, V_t)$. It is assumed to be absolutely continuous, i.e., if

$$P_t^i(k \times [0, y]) := P^i\{Z_t = k, V_t \in [0, y]\}, \tag{3.9}$$

then

$$\pi_k^i(t, u) du = P\{Z_t = k, V_t \in [u, u + du)\}, \; i, k \in E, u, t \in \mathbb{R}_+. \tag{3.10}$$

Since all pertinent processes in equilibrium have been studied, the initial state is not considered. So the index $(i)$ from $\pi$ is dropped. To find the limiting distribution $\boldsymbol{\pi} = \{\pi_0, \pi_1, \ldots, \pi_w\}$ of the process $(Z_t)$, we will start with the Kolmogorov differential equations, assuming that $A(x)$ has the Radon-Nikodym density $a(x)$, which, in addition, is pointless continuous. From (3.10), the Kolmogorov equations are as follows:

$$\left(\frac{\partial}{\partial t} - \frac{\partial}{\partial u}\right)\pi_k(t,u) = \begin{cases} -\lambda\pi_k(t,u) + \lambda\pi_{k-1}(t,u) + \pi_{k-1+v}(t,0)a(u), \\ \qquad k = 1,\ldots,w-v, \\ -\lambda\pi_k(t,u) + \lambda\pi_{k-1}(t,0) + \pi_w(t,0)a(u), \\ \qquad k = w-v+1,\ldots,w-1, \\ -\lambda\pi_w(t,u) + \lambda\pi_{w-1}(t,0), \\ \qquad k = w. \end{cases} \quad (3.11)$$

The process $(Z_t, V_t)$ is also semi-regenerative relative to the sequence $T_0, T_1, \ldots$ and its limiting distribution exists [7-8]. Let

$$\pi_k(u) := \lim_{t \to \infty} \pi_k(t, u), \tag{3.12}$$

$$\widehat{\pi}_k(\theta) := \int_{\mathbb{R}_+} \pi_k(u) e^{-\theta u} du, \ \mathrm{Re}(\theta) \geq 0, \ k \in E. \tag{3.13}$$

Letting $t \to \infty$ in (3.11) and then applying the Laplace transform we have:

$$(\theta - \lambda)\widehat{\pi}_k(\theta) = \pi_k(0) - \lambda\widehat{\pi}_{k-1}(\theta) - \pi_{k+v-1}(0)\alpha(\theta), \tag{3.14}$$
$$k = 1, \ldots, w-v-1;$$
$$(\theta - \lambda)\widehat{\pi}_k(\theta) = \pi_k(0) - \lambda\widehat{\pi}_{k-1}(\theta) + \pi_w(0)\alpha(\theta), \tag{3.15}$$
$$k = w-v, \ldots, w-1;$$
$$(\theta - \lambda)\widehat{\pi}_w(\theta) = \pi_w(0) - \lambda\widehat{\pi}_{w-1}(\theta). \tag{3.16}$$

With $\theta \downarrow 0$ in (3.14)-(3.16),

$$-\lambda\pi_k = \pi_k(0) - \lambda\pi_{k-1} - \pi_{k+v-1}(0), \ k = 1, \ldots, w-v-1, \tag{3.17}$$

$$-\lambda\pi_k = \pi_k(0) - \lambda\pi_{k-1} - \pi_w(0), \ k = w-v, \ldots, w-1, \tag{3.18}$$

$$-\lambda\pi_w = \pi_w(0) - \lambda\pi_{w-1}. \tag{3.19}$$

Once, $\pi_k(0)$, $k = 0, 1, \ldots, s$ are revealed, $\pi_k$, $k = 0, \ldots, w$ could be found accordingly after calculating the invariant probability measure $\mathbf{P} = (P_0, P_1, \ldots, P_w)$ of $\{X_n\}$ from the TPM of the system which covered on the next section. To find the unknowns $\pi_k(0)$ when $k = 0, 1, \ldots, w$. Denote

$$p_{jk}(t) := P^j\{Z_t = k \,|\, T_1 > t\}. \tag{3.20}$$

Note that

$$\int_{\mathbb{R}_+} p_{jk}(t)\, A(dt) = p_{jk}, \ i, k \in E, \tag{3.21}$$

are the transition probabilities of the embedded MC $(X_n)$. We use the natural assumption that

$$P^j\{Z_t = k \,|\, T_1 > t\} = P^j\{Z_t = k \,|\, T_1 > t + y\}, \ \forall y \geq 0.$$

Let $K_t^j(k \times [0, y]) := P^j\{Z_t = k, V_t \in [0, y], T_1 > t\}$ and from the main convergence theorem,

$$\lim_{t\to\infty} P^i_t(k \times [0,y]) = \frac{1}{a}\sum_{j\in E} P_j \int_{\mathbb{R}_+} K^j_t(k \times [0,y])dt, \quad k \in E, \, y \in \mathbb{R}_+. \tag{3.22}$$

By elementary probability arguments,

$$K^j_t(k \times [0,y]) = p_{jk}(t)[A(t+y) - A(t)]$$

and then

$$\lim_{t\to\infty} P^i_t(k \times [0,y]) = \lim_{t\to\infty} P^i\{Z_t = k, V_t \in [0,y]\}$$

$$= \frac{1}{a}\int_{x=0}^{y} \sum_{j\in E} P_j \int_{\mathbb{R}_+} p_{jk}(t) a(t+x) \, dt \, dx. \tag{3.23}$$

From (3.23),

$$\pi_k(x) = \frac{1}{a}\sum_{j\in E} P_j \int_{\mathbb{R}_+} p_{jk}(t) a(t+x) \, dt. \tag{3.24}$$

When $x \to 0$ in (3.24), due to $(\mathbf{P}, \mathbf{1}) = 1$ and (3.7),

$$\pi_k(0) = \frac{P_k}{a}. \tag{3.25}$$

### 3.1. Process of Model 2 Type 1
The Type 1 of Model 2 is the case that the total number of services $w$ in the system is larger than double of the batch independent supplier size $v$ (i.e., $w \geq 2v$). In the Type 1, The transition probability are determined from (3.2) and the TPM is constructed as follows:

$$\mathbf{M}_{(w+1;\, w\geq 2v)} = \begin{pmatrix} \psi(0) & \psi(1) & \cdots & \cdots & \psi(w-v) & 0 & \cdots & 0 \\ \vdots & \vdots & \vdots & \vdots & \vdots & \vdots & \vdots & \vdots \\ \psi(0) & \psi(1) & \cdots & \cdots & \psi(w-v) & 0 & \cdots & 0 \\ 0 & \psi(0) & \vdots & \cdots & \psi(w-1-v) & 0 & \cdots & 0 \\ \vdots & \vdots & \vdots & \vdots & \vdots & 0 & \cdots & 0 \\ 0 & 0 & \vdots & \cdots & \vdots & \vdots & \vdots & \vdots \\ 0 & 0 & \psi(0) & \cdots & \psi(w-2v) & \cdots & \cdots & 0 \\ 0 & \cdots & \cdots & 0 & 0 & 0 & \cdots & 0 \\ 0 & 0 & \cdots & \cdots & 0 & 0 & \cdots & 0 \end{pmatrix} \tag{3.26}$$

and the invariant probability measure $\mathbf{P} = (P_0, P_1, \ldots, P_w)$ of $\{X_n\}$ is determined from $(\mathbf{P}, \mathbf{1}) = 1$ and (3.7)-(3.8). From (3.17), summation over $j = 1, \ldots, k \, (\leq w - v)$ in (3.17):

$$\lambda(\pi_k - \pi_0) = \begin{cases} \displaystyle\sum_{j=1}^{k} \pi_{v+j-1}(0) - \widehat{\mathcal{B}}_k, & k = 1, \ldots v, \\ \displaystyle\sum_{j=v+1}^{k-v} \pi_{v+j}(0) - \widehat{\mathcal{B}}_{v-1}, & k = v+1, \ldots, w-v, \end{cases} \tag{3.27}$$

where

$$\widehat{\mathcal{B}}_l = \sum_{i=1}^{l} \pi_i(0) = \frac{1}{a}\sum_{i=1}^{l} P_i \qquad (3.28)$$

Summing up the equation of (3.18)-(3.19) in $j$ from $w-v+1$ to $s$ then yields:

$$\lambda(\pi_k - \pi_0) = (k-(w-v))\pi_w(0) - \widehat{\mathcal{B}}_{v-1} \qquad (3.29)$$

The continuous time parameter process of Type 1 (Model 2) could be analyzed from (3.25) and (3.27)-(3.29) and limiting distribution $\pi$ is as follows:

$$\pi_n = \mathcal{G}_n + \pi_0, \ n = 1, \ldots, w \qquad (3.30)$$

and using $(\pi, \mathbf{1}) = 1$,

$$\pi_0 = 1 - \sum_{n=1}^{w} \pi_n, \qquad (3.31)$$

where

$$\lambda \cdot \mathcal{G}_n = \begin{cases} \frac{1}{a}\left(\sum_{j=1}^{n} P_{v+j-1}\right) - \widehat{\mathcal{B}}_n, & n = 1, \ldots, v, \\ \frac{1}{a}\left(\sum_{j=v+1}^{n-v} P_{v+j}\right) - \widehat{\mathcal{B}}_{v-1}, & n = v+1, \ldots, w-v, \\ \frac{1}{a}(n-w+v)P_w - \widehat{\mathcal{B}}_{v-1}, & n = w-v+1, \ldots, w. \end{cases} \qquad (3.32)$$

From (3.26), (3.29) and (3.30)-(3.32),

$$\pi_0 = 1 - \left(\sum_{k=1}^{w} \mathcal{G}_k - w \cdot \pi_0\right),$$

therefore,

$$\pi_0 = \left[1 - \sum_{k=1}^{w} \mathcal{G}_k\right]/(1+w). \qquad (3.33)$$

For the process $Z_t^1$, the corresponding formulas yield for the limiting distribution $\pi^1$ $\left(\pi_k^1 = \lim_{t\to\infty} P\{Z_t^1 = k\}, \ k = 0, \ldots, w\right)$ is as follows:

$$\pi_k^1 = \pi_{w-k}, k = 0, 1, \ldots, w, \qquad (3.34)$$

along with (2.6).

### 3.2. Process of Model 2 Type 2
Alternatively, the embedded process $\{X_n\}$ of Model 2 is changed when $w < 2v$ (Model 2 Type 2) and the TPM becomes as follows from (3.2a) and (3.3):

$$\mathbf{M}_{(w+1;\, w<2v)} = \begin{pmatrix} \psi(0) & \psi(1) & \ldots & \psi(w-v) & 0 & \ldots & 0 \\ \vdots & \vdots & \vdots & \vdots & \vdots & \vdots & \vdots \\ \psi(0) & \psi(1) & \ldots & \psi(w-v) & 0 & \ldots & 0 \\ 0 & 0 & \ldots & 0 & 0 & \ldots & 0 \\ \vdots & \vdots & \vdots & \vdots & \vdots & \vdots & \vdots \\ 0 & 0 & 0 & \ldots & 0 & \ldots & 0 \end{pmatrix}, \qquad (3.35)$$

and the TPMs from (3.35) repent the same finite sub-matrix which occurs within the corresponding TPM of $M/G^v/1/\infty$ queueing system (3.4). Similarly, the continuous time parameter process for Type-2 could be found from (3.17). From (3.17)-(3.20), summation over $j = 1, \ldots, k\ (\leq v)$ in (3.18)-(3.19):

$$\lambda(\pi_k - \pi_0) = \begin{cases} \sum_{j=1}^{k} \pi_{v+j-1}(0) - \widehat{\mathcal{B}}_k, & k = 1, \ldots w-v, \\ \sum_{j=1}^{w-v} \pi_{v+j-1}(0) + (k - (w-v))\pi_w(0) - \widehat{\mathcal{B}}_k, & k = w-v+1, \ldots, w, \end{cases} \qquad (3.36)$$

and summing up the equation of (3.18)-(3.19) in $j$ from $v+1$ to $w$ then yields:

$$\lambda(\pi_k - \pi_0) = \pi_{v-1}(0) + \sum_{j=k-v+1}^{w-v} \pi_{v+j}(0) + (k - (w-v))\pi_w(0) - \widehat{\mathcal{B}}_{v-1}. \qquad (3.37)$$

From (3.25), (3.28) and (3.30)-(3.31), we finally arrive at the limiting distribution $\boldsymbol{\pi}$ with the different $\mathcal{G}_n$:

$$(3.38)$$

$$\lambda \cdot \mathcal{G}_n = \begin{cases} \frac{1}{a}\left(\sum_{j=1}^{n} P_{v+j-1}\right) - \widehat{\mathcal{B}}_n, & n = 1, \ldots, w-v, \\ \frac{1}{a}\left(\sum_{j=1}^{w-v} P_{v+j}\right) + \frac{1}{a}(n-w+v)P_w - \widehat{\mathcal{B}}_n, & n = w-v+1, \ldots, v, \\ \frac{1}{a}\left(\sum_{j=n-v+1}^{w-v} P_{v+j}\right) + \frac{1}{a}(n-w+v)P_w - \widehat{\mathcal{B}}_{v-1}, & n = v+1, \ldots, w. \end{cases}$$

For the process $Z_t^1$ for Type 2, the corresponding formulas yield same from (3.30) and (3.33)-(3.34) but with the different parameters from (3.38).

## 4. THE OPTIMALITY OF THE SHARING ECONOMY SYSTEM

In this section, a class of optimization problem that arise a stochastic sharing business operations is considered. Let us formalize a pertinent optimization problem. Let a strategy, say $\Psi$, specify, ahead of the time, a set of acts that imposed on the queueing system, such as the choice of the batch size of an individual contracts ($v$) in the service platform, the total number of resource suppliers ($w$) and so on. Denote by $\phi(\Psi, C, t)$ the expected costs within $[0, t]$, due to the strategy $\Psi$ and costs $C$ and define the expected cumulative cost rate over an infinite horizon:

$$\phi(\Psi, C) := \lim_{t \to \infty} \frac{1}{t} \phi(\Psi, C, t). \tag{4.1}$$

Let $h(n)$ denote the cost function associated with the holding the contract in the (service) platform for $n$ contractors which are remained in the supplier pool. If $h(u)$, a linear function, i.e., $h(n) = c_H \cdot n$, then the expected holding cost for contractors in the supplier pool during the interval $[0, t]$ is

$$Uh(t) = \mathbb{E}\left[h\left(\int_0^t h(Z_u^1)du\right)\right] = \mathbb{E}\left[h(k) \int_{u=0}^t \mathbf{1}_{\{k\}}(Z_u^1)du\right] \tag{4.2}$$

(which by Fubini's Theorem is)

$$= c_H \cdot \left(\sum_{k=0}^s \int_{u=0}^t k P\{Z_u^1 = k\}du\right). \tag{4.2a}$$

Since $Z_t^1$ is the number of available suppliers on the service platform, $M_t = (Z_t^1 - v)^+$ gives the number of suppliers which are reserved by the company (i.e., platform provider) for urgent cases (i.e., individual contractors are not enough for support all consumers). The expected cost for reserved suppliers in the interval $[0, t]$ is

$$Ug(t) = \mathbb{E}\left[\int_0^t g(M_u)du\right] = \mathbb{E}\left[g((k-v)^+) \int_{u=0}^t \mathbf{1}_{\{k\}}(Z_u^1)du\right] \tag{4.3}$$

(which by Fubini's Theorem is)

$$= \sum_{k=v+1}^w g(k-v) \int_{s=0}^t P\{Z_s^1 = k\}ds,$$

if $g(n)$ is a linear function (i.e., $g(n) = c_R \cdot n$),

$$= c_R \cdot \left(\sum_{k=v+1}^w \int_{u=0}^t (k-v) P\{Z_u^1 = k\}du\right). \tag{4.3a}$$

Since the Markov renewal function $D(t) = \mathbb{E}\left[\sum_{n \geq 0} \mathbf{1}_{[0,t]}(T_n)\right]$ represents the total number of services which consumed by customers in the time interval $[0, t]$ where $\mathbf{1}_B$ is the indicator function of a set $B$. The functional $Ud(u) = c_D(\lambda/v)D(t)$ gives the set cost for posting the available suppliers on the service platform. The cumulative cost of the entire procedures is involved in the sharing service operations in the interval $[0, t]$ is

$$\begin{aligned}\phi(\Psi, C, t) &= Uh(t) + Ug(t) + Ud(t) \\ &= \left(\sum_{k=0}^w \int_{u=0}^t h(k) P\{Z_u^1 = k\}du\right) \\ &+ \left(\sum_{k=v+1}^w \int_{u=0}^t g(k-v) P\{Z_u^1 = k\}du\right) + c_D(\lambda/v)D(t).\end{aligned} \tag{4.4}$$

Now we turn to convergence theorems for Markov renewal and semi-Markov processes [10],

(i) $\lim_{t\to\infty} \frac{1}{t} D(t) = \frac{1}{a}$

(ii) $\lim_{t\to\infty} \frac{1}{t} \int_0^t P^i\{Z_s^1 = k\} ds = \pi_k^1$

to arrive at the objective function $\phi(\Psi, C)$, which gives the total expected rate of all processes over an infinite horizon. In light of equations (i-ii) we have

$$\phi(\Psi, C) = \lim_{t\to\infty} \frac{1}{t} \phi(\Psi, C, t) = \sum_{k=0}^{w} h(k)\pi_k^1 + \sum_{k=v+1}^{w} g(k-v)\pi_k^1 + c_D\left(\frac{\lambda}{v}\right) D(t). \quad (4.5)$$

With the cost functions being linear functions, we have

$$\phi(\Psi, C) = c_H \cdot \mathbb{E}[Z_\infty^1] + c_R \left[\sum_{k=v+1}^{w} (k-v)\pi_k^1\right] + c_D\left(\frac{\lambda}{v}\right)\left[\frac{1}{a}\right]. \quad (4.6)$$

Recall that from (3.34) we have

$$\begin{cases} \pi_n^1 = \mathcal{G}_{w-n} + \pi_w^1, & n = 0, 1, \ldots, w-1, \\ \pi_w^1 = \left[1 - \sum_{k=1}^{w} \mathcal{G}_k\right]/(1+w). \end{cases} \quad (4.7)$$

With these attachments and (3.4) and (3.7a), the formula (4.6) for the objective function is complete.

## 5. SUPPLIER CONTRACT OPTIMIZATION CASE

On a special model of the sharing business operation situation, this section intends to demonstrate the tractability of our results in sections 2-4. We allow an exponential distribution of inter-arrival duration between (individual) contract posting on the platform. Since the method involves operating with an embedded Markov process in single channel open queue with a finite buffer, we get back to section 3 for some particular formulas.

### 5.1. $M^v/M/1/w$ Sharing Business Operations

The pertinent special formulas of section 4 under these assumptions are as follows:

$$A(x) = 1 - e^{-\frac{1}{a}x}, \quad (5.1)$$

$$\alpha(\theta) = \int_{x=0}^{\infty} \left(\frac{1}{a}\right) e^{-(\theta+\frac{1}{a})x} dx = \frac{1}{1+a\theta}, \quad (5.2)$$

$$\alpha(\lambda(1-z)) = \frac{1}{1+\lambda a(1-z)}, \quad (5.3)$$

$$Q(z) = \sum_{i\geq 0} Q_i \cdot z^i = \frac{\left[\sum_{i=0}^{v-1} Q_i(z^v - z^i)\right]}{(1-\lambda a)z^v - \lambda a z^{v+1} - 1}. \quad (5.4)$$

From Chaudhry and Templeton [6],

$$Q_i = \left(1 - \frac{1}{z_0}\right)\left(\frac{1}{z_0}\right)^i \tag{5.5}$$

where, from (5.4),

$$z_0 \in \{\exists z : \{(1-\lambda a)z^v - \lambda a z^{v+1} - 1 = 0\} \wedge \{z > 1\}\}.$$

The remaining parameter such as $\widehat{\mathcal{B}}_l$ from (3.28), $\mathcal{G}_k$ from (3.32) (or (3.38)) and $P_k$ from (3.7) to get $\pi_k$ from (3.30)-(3.31) could be calculated accordingly. After finding $\pi_k$ ($k = 0, \ldots, w$), we can find $\pi_k^1$, $k = 0, \ldots, w$ from (4.7).

## 5.2 Optimization with the two types of supplier contracts

For our optimization example, we specify the remaining three primary cost functions are as follows:

$$h(n) = c_H \cdot n, \; g(n) = c_R \cdot n, \tag{5.6}$$

where $c_H$ is the holding cost per contract in the service platform, $c_R$ is the supplier costs which reserved by the company for urgent cases. From (4.5) and (5.6), we have

$$\phi(\Psi, C) = c_H \cdot \left[\overline{Z}_\infty^1\right] + c_R \cdot \left[\sum_{k=v+1}^{w}(k-v)\pi_k^1\right] + c_D\left(\frac{\lambda}{v}\right)\left[\frac{1}{a}\right], \tag{5.7}$$

where

$$\overline{Z}_\infty^1 := \mathbb{E}[Z_\infty^1] = \sum_{k=0}^{w} k\pi_k^1.$$

Finally, we arrive at the following expression for the objective function:

$$\phi(\Psi(v), C) = c_H \cdot \left[\overline{Z}_\infty^1\right] + c_R \cdot \left[\sum_{k=v+1}^{w}(k-v)\pi_k^1\right] + c_D\left(\frac{\lambda}{v}\right)\left[\frac{1}{a}\right]. \tag{5.8}$$

Here, we use formulas (4.6)-(4.7) and (3.32) for Type-1 (or (3.38) for Type-2). Notice that we have the parameter $v$ vary. We restrict the initial strategy of this model to one, which includes only the control level $v$ of the batch size of each arrival. In other words, we need to find a $v$ such that

$$\phi(\Psi(v), C) = min\{\phi(\Sigma(v), C) : b = 1, 2, \ldots, v_{max} \; (< w)\}. \tag{5.9}$$

As an illustration in Figure 2, let us take $c_H = 3$, $c_R = 1$ and $c_D = 80$. The inter-duration between posting one bulk of contracts on the (online) platform is exponentially distributed with the mean $a = 1.3$ and the parameter of consuming services by customers is $\lambda = 2.2$. Take the maximum number of contracts (suppliers) is 35 and the maximal bulk size of the individual contracts must be smaller than the size of the total capacity $w$ (i.e., $v_{max} \leq w$). Now, we calculate $\phi(\Psi(v), C)$ and $v_0$ that gives a minimum for $\phi(\Psi(v), C)$. In other words, the control level $v_0$ stands for the optimal batch size of the individual supplier contracts which minimizes the total operating cost of the sharing business. Below is a plot of $\phi(\Psi(v), C)$ for $N \in \{1, 2, \ldots, w(= 35)\}$.

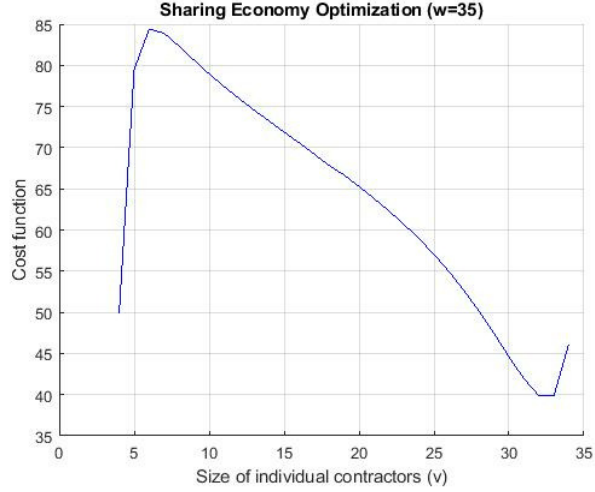

**Figure 2.** Optimization Example

Our calculation yields that $v_0 = 33$ for which the minimal cost equals 39.8673. When the total number of the allowed suppliers is 35, we obtain the threshold value $v_0 = 33$ which gives us the decision point that is the number of bulk contracts for the individual suppliers to minimize the cost of the sharing economy system operations. It also indicates that only 2 suppliers are reserved by the platform provider (i.e., $s_0 = w - v_0$). In addition, it is also feasible to visualize the cost objective function with two variables in the platform (see Figure 3).

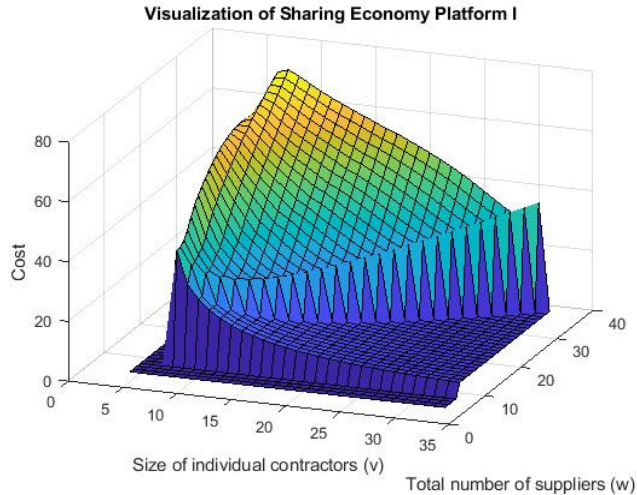

**Figure 3.** Visualization of the cost analysis based on two variables

As a part of reasonable performance measures, let us consider a capability factor $\varrho$, which represents the estimated loss rate of the potential service (customers) at any moment of time in equilibrium:

$$\varrho = \left[\frac{\lambda a}{w} - 1\right]^+. \qquad (5.10)$$

In this particular case, this value is zero (i. e., $\varrho = 0$) and it indicates that no customer losses in the system. This value could not only be a performance measure but be also a constraint to optimize an objective function of the sharing economy platform.

## 6. CONCLUSIONS

The paper has propounded a compounding model of operating the sharing economy system by using the $G^v/M/1$ queueing model with the stochastic congruence. For each of these systems, we have determined the stationary distributions for the embedded and the continuous time parameter processes. The queueing process in the sharing economy platform coincides with a conventional $M/G^v/1$ queueing system with a finite waiting room. The theoretical approach for a sharing economy platform could become the basis for resolving social and economic challenges by better understanding before adapting complex intelligent systems. Analytically tractable results are obtained by using a stochastic congruence, semi-regenerative analysis and semi-Markov process. This analytic approach supports the theoretical background of operating the sharing economy service platform. The comparison and the gap analysis between the theoretical solutions and the results based on the actual data from real sharing economy firms might be considered as a future research topic in this area of studies. Additionally, the theoretical mathematical modeling for various social challenges could be yet another research topics in the near future.

## ACKNOWLEDGMENT

The authors are most thankful to the referees whose comments are very constructive. Additionally, the authors acknowledge support from the Center for Cyber-Physical Systems, Khalifa University, under Grant Number 8474000137-RC1-C2PS-T3 and there is no available data to be stated.